\documentclass[a4paper,11pt,reqno]{amsart}
\usepackage[margin=2.7cm]{geometry}
\usepackage{cite}
\usepackage{newpxtext}
\usepackage{cabin}
\usepackage[varqu,varl]{inconsolata}
\usepackage[bigdelims,vvarbb]{newpxmath}
\usepackage[cal=cm,bb=esstix,bbscaled=1.05]{mathalfa}
\usepackage{booktabs}
\usepackage{caption}
\usepackage{enumitem}
\usepackage[colorlinks]{hyperref}
\usepackage{xcolor}
\usepackage{graphicx}
\hypersetup{
    colorlinks,
    linkcolor={red!50!black},
    citecolor={blue!50!black},
    urlcolor={blue!80!black}
}

\newcommand{\orcid}[1]{\,\resizebox{8px}{!}{\href{https://orcid.org/#1}{\includegraphics{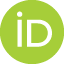}}}}
\newcommand{\PP}{\mathbb P}

\newcommand{\Q}{\mathbb Q}
\newcommand{\C}{\mathbb C}
\newcommand{\R}{\mathbb R}
\newcommand{\Z}{\mathbb Z}
\newcommand{\cO}{\mathcal O}
\newcommand{\tf}{3\nobreakdash-fold}
\newcommand{\tfs}{3\nobreakdash-folds}
\newcommand{\ff}{4\nobreakdash-fold}
\newcommand{\ffs}{4\nobreakdash-folds}
\newcommand{\nf}{$n$\nobreakdash-fold}
\newcommand{\nfs}{$n$\nobreakdash-folds}
\renewcommand{\phi}{\varphi}
\newcommand{\eps}{\varepsilon}
\DeclareMathOperator{\GL}{GL}
\DeclareMathOperator{\Aut}{Aut}
\DeclareMathOperator{\shed}{shed}
\DeclareMathOperator{\vertices}{vertices}
\newcommand{\De}{\Delta}
\newcommand{\si}{\sigma}
\newcommand{\Si}{\Sigma}
\newcommand{\Ga}{\Gamma}
\newcommand{\Qfact}{$\Q$\nobreakdash-factorial}
\numberwithin{equation}{section}

\begin{document}
\title{Toric Sarkisov links of toric Fano varieties}
\author[G.\ Brown]{Gavin Brown\orcid{0000-0002-4087-5624}}
\address{Mathematics Institute\\Zeeman Building\\University of Warwick\\Coventry\\CV4 7AL\\UK}
\email{G.Brown@warwick.ac.uk}
\author[J.\ Buczy\'nski]{Jaros{\l}aw Buczy\'nski\orcid{0000-0002-9811-4238}}
\address{Instytut Matematyczny PAN\\ul.\ \'Sniadeckich 8\\00-656 Warszawa\\Poland}
\email{jabu@mimuw.edu.pl }
\author[A.\ M.\ Kasprzyk]{Alexander Kasprzyk\orcid{0000-0003-2340-5257}}
\address{School of Mathematical Sciences\\University of Nottingham\\Nottingham\\NG7 2RD\\UK}
\email{a.m.kasprzyk@nottingham.ac.uk}
\begin{abstract}
We explain a web of Sarkisov links that overlies the classification of Fano weighted projective spaces in dimensions~$3$ and~$4$, extending results of Prokhorov.
\end{abstract}
\maketitle
\section{Introduction}
A normal projective~$n$-dimensional complex variety is called~\emph{Fano} if it has ample anticanonical class and canonical singularities. The construction and classification of Fano \nfs\ is a major concern of birational geometry. Dimensions~$n=3$ and~$4$ are at the cutting edge, where birational methods of construction play a central role. In various respects, toric Fano varieties arise at the extremes of classification (compare~\cite{grdb}). In this paper we consider the toric birational geometry of certain toric Fano \ffs, and in particular the Sarkisov links between them; we review terminology in~\S\ref{sec!sarkisov} and the results are surveyed in~\S\ref{sec!dim4}, with full details of over a million links relegated to a webpage~\cite{grdbweb}. Although as a study of the 4-dimensional Sarkisov program this is a baby case, it does provide a large number of examples and gives some indication of the \ff\ phenomena we can expect to encounter when we push beyond toric. It also describes classes of new toric Fano \ffs, as we discuss in~\S\ref{sec!midpoints}.

A~\emph{Fano polytope} is a polytope~$\De\subset \R^n=L\otimes\R$ whose vertices lie on a fixed lattice~$L\cong\Z^n$ and whose only strictly interior lattice point~$p\in\De^\circ\cap L$ is the origin. A Fano polytope is~\emph{terminal} if its vertices are the only lattice points lying on its boundary, and it is~\emph{\Qfact} (or~\emph{simplicial}) if each facet is an~$(n-1)$-simplex. In~$3$~dimensions Fano polytopes are classified up to~$\Aut(L)$ into~$674,\!688$ cases~\cite{tf3,kasprzyk10}:~$12,\!190$~are simplicial;~$634$~are terminal;~$233$~are~\Qfact\ terminal, and of these only~$8$ have the minimal number of four vertices, the terminal simplicial tetrahedra. In~$4$~dimensions, Fano polytopes are not yet classified, but the terminal simplicial polytopes (on five vertices) are, with~$35,\!947$ cases~\cite{kasprzyk13}.

Toric geometry connects these two realms. Given a Fano polytope~$\De$, the spanning fan~$\Si$ of~$\De$ is the collection of lattice cones on the facets and smaller strata of~$\De$. This fan~$\Si$ determines a normal complex projective variety~$X=X_\Si$~\cite{danilov78,fulton}, that has ample anticanonical class~$-K_X$ (by the convexity of~$\De$) and canonical singularities (by the assumption on~$\De^\circ\cap L$)~\cite{reid82}. In other words,~$X$ is a Fano \nf\ with canonical singularities, and furthermore it is terminal or~\Qfact\ exactly when~$\De$ is. In the \Qfact\ case, the Picard rank~$\rho_X$ of~$X$ is the number of vertices of~$\De$ minus~$n$. A~\emph{Mori\nobreakdash--Fano \nf} is a \Qfact\ terminal Fano \nf\ $X$ with~$\rho_X=1$. The toric Mori\nobreakdash--Fano \nfs\ are those weighted projective spaces and their quotients by the finite discriminant group~$L/\left<\vertices(\De)\right>$, often called~\emph{fake weighted projective spaces} (compare~\cite[6.4]{wero}), that have terminal singularities. Thus there are precisely~$8$ toric Mori\nobreakdash--Fano \tfs\ (compare~\cite{bb92}), and~$35,\!947$ toric Mori\nobreakdash--Fano \ffs\ (of which~$24,\!511$ are weighted projective spaces).

Toric maps from~$X$ are controlled by the combinatorics of the fan~$\Si$. Prokhorov~\cite[\S10]{prokhorov10} works out a beautiful web of Sarkisov links between~$7$ of the weighted projective space \tfs. In this paper we complete that web with all possible toric links. We then extend the programme into dimension~$4$, closely following~\cite{reid82}, giving examples to illustrate the phenomena that arise, and providing full results online at~\cite{grdbweb} (with associated Magma code~\cite{magma} that generates them). Guerreiro~\cite{tiago} recently computed all Sarkisov links from~$\PP^4$ that start with a weighted point blowup, which we discuss in~\S\ref{sec!highdiscrep}.

Before starting, we return to the polytope point of view.
Reid~\cite[(4.2)]{reid82} defines the~\emph{shed} of a fan~$\Si$ to be
\[
\shed\Si = \!\!\bigcup_{C\in\Si^{(n)}} \!\!\shed C
\]
where the union is over the top-dimensional cones~$C$ of~$\Si$, and~$\shed C$ is the convex hull of the primitive points on each (extreme) ray of~$C$ together with the origin. When~$\Si$ is the spanning fan of a Fano polytope~$\De$, then~$\shed\Si=\De$; in contrast, if~$\Si$ is the fan of~$X$ with~$-K_X$ not nef, then the shed is not convex; see~\cite[(4.3)]{reid82} or the shed of~$Y_2$ in Figure~\ref{fig!pic} below. In these terms, the operations~(\ref{sec!operations}.\ref{type:extraction}--\ref{type:divisor}) below are lattice piecewise linear~(LPL) surgeries on the shed, adding or removing vertices, or snapping~LPL portions of the boundary into other~LPL configurations, all the time preserving the property that the only lattice points in the shed are its vertices and the internal origin. A Sarkisov link then appears as a kind of~LPL origami, moving from one Fano polytope to another by a series of~LPL modifications, as we illustrate in Figure~\ref{fig!pic} below.
\section{The context}
\subsection{Overview of the Sarkisov program}\label{sec!sarkisov}
A morphism~$\phi\colon V\rightarrow S$ of normal projective varieties is, by definition, a Mori fibre space if~$\phi_*\cO_V=\cO_S$,~$V$ has \Qfact\ terminal singularities,~$\dim S < \dim V$,~$\rho_V=\rho_S+1$ and~$-K_V$ is~$\phi$\nobreakdash-ample. The Sarkisov program~\cite{corti95,HMsarkisov} decomposes birational maps~$f\colon V \dashrightarrow V'$ between Mori fibre spaces. It factorises any such birational map into a composition of very particular birational maps called~\emph{Sarkisov links}, which we explain in our context below, and so understanding individual Sarkisov links becomes crucial. Sarkisov links are built by patching together steps of the form
\begin{equation}\label{eq!square}
\begin{array}{ccccc}
&&U\\
&\swarrow&&\searrow\\
V_1&&&&V_2\\
&\searrow&&\swarrow\\
&&S
\end{array}
\end{equation}
The key point is that~$\rho_{U/S}=2$ in the square~\eqref{eq!square}, so that up to isomorphism~$U$ admits at most two~$S$-morphisms of relative Picard rank~$1$ with connected fibres. The construction of a Sarkisov link is an inductive procedure in which the map~$U\rightarrow V_1$ is given and we must solve for~$U\rightarrow V_2$ (over~$S$). The link is called~\emph{bad}, and we abandon it, if~$U\rightarrow V_2$ does not exist, or if the properties of~$V_2$ break any running hypotheses, for example on the singularities permitted. (See~\cite[2.2]{corti00} for full details in these terms.)

When~$X$ is a Mori\nobreakdash--Fano \nf, a Sarkisov link will always take the form
\begin{equation}\label{eq!link}
\begin{array}{cccccccccccc}
&&Y_1 & \dashrightarrow & Y_2 \dashrightarrow & \cdots & \dashrightarrow & Y_s \\
& \swarrow && \searrow\quad\swarrow &&&&& \searrow \\
X &&& Z_1 &&&&&& Z_s = X'
\end{array}
\end{equation}
where~$g_1\colon Y_1 \rightarrow X$ is an extremal extraction (that is,~${g_1}_*\cO_{Y_1}=\cO_X$,~$Y_1$ has terminal \Qfact\ singularities,~$\rho_{Y_1/X}=1$ and~$-K_{Y_1}$ is relatively ample), each~$Y_i\dashrightarrow Y_{i+1}$ is a generalised flip (that is, a flip or a flop or an antiflip) that we may factorise as a pair of birational morphisms $Y_i\rightarrow Z_i\leftarrow Y_{i+1}$. When~$Y_i\dashrightarrow Y_{i+1}$ is an antiflip, the link is bad if~$Y_{i+1}$ does not have terminal singularities. The final extremal contraction~$Y_s\rightarrow X'$ may be to another Fano \nf~$X'$, or have~$\dim Z_s<\dim Y_s$ so that~$Y_s\rightarrow Z_s$ is a strict Mori fibre space; the link is referred to as being of Type~I or Type~II in these two cases respectively~\cite{corti00} (the opposite naming convention is used in~\cite{HMsarkisov}).
\subsection{Hard Fano varieties from easy ones}
Sarkisov links~\eqref{eq!link} of Type~I relate certain members of one family of Fano varieties with certain members of another. In~\S\ref{sec!dim3}, we calculate all Sarkisov links composed of toric maps between any two toric Fano \tfs. The results are listed in Table~\ref{fig!links}, and our choice of presentation illustrates how we might treat this idea as passing from simpler varieties to more complicated ones. Corti--Mella~\cite{CM} describe links from (general) quasismooth codimension~$2$ Fano \tf\ to a (particular) singular Fano hypersurface. Those links start with a projection. Reversing that step as an `unprojection' is the key idea behind the constructions in codimension~$4$ of~\cite{BKR}, which are extended to Sarkisov links in~\cite{Campo}. One may view this construction and classification method, and results of Takeuchi, Takagi and others~\cite{takeuchi,takagi}, in this same light: crudely speaking, they list the possible numerical data that could be realised by a Sarkisov link, and then construct or eliminate each case.

In three dimensions, the Graded Ring Database~\cite{grdbweb} presents a first overview of the numerical data of the classification of Fano \tfs, again, crudely speaking, listing the possible numerical data that a Fano \tf\ might have, with the aim of constructing or eliminating each case; see~\cite[Fig.~1]{grdb}. In that picture, families are connected by projection, but the numerical data of Sarkisov links, such as~\cite[(2.8)]{takeuchi} or~\cite[0.3]{takagi}, is missing. We regard our exercise here is a first (very baby) step in understanding how the possible numerical data of Sarkisov links between Fano varieties might be included in the Graded Ring Database.

\begin{table}[tbp]
\[
\small
\begin{array}{c|c|c|ccc|c}
\toprule
\# & X& \text{blowup} & \text{antiflip} & \text{flop / flip} & \text{blowdown} & \text{end of link} \\
\cmidrule(lr){1-7}
1 &\PP(1,1,1,1) & (1,1,0) &&&\text{Mfs}& \PP^2 / \PP^1 \\
2 && (1,1,1) &&&\text{Mfs}& \PP^1 / \PP^2 \\
3 && (1,1,2) &&& (1,1,0) & \PP(1,1,1,2) \\
4 && (1,2,3) & 1,1,-1,-2 && (1,1,2) & \PP(1,1,2,3) \\
5 && (1,2,5) & 1,1,-1,-4 && \tfrac{1}3(1,1,2) & \PP(1,3,4,5) \\
\cmidrule(lr){1-7}
6 & \PP(1,1,1,2) & \tfrac{1}2(1,1,1) &&&& \PP^1 / \PP^2 \\
7 && (1,1,1) && 2,1,-1,-1 & \text{Mfs} & \PP^2 / \PP^1 \\
8 && (1,1,2) & 2,1,-1,-3 && (1,1,1) & \PP(1,1,2,3) \\
9 && (1,1,2) &&&\text{Mfs}& \PP^1 / \PP(1,1,2) \\
10 && (1,1,3) &&& (1,1,0) & \PP(1,1,2,3) \\
11 && (1,1,3) & 2,1,-1,-5 && \tfrac{1}2(1,1,1) & \PP(1,2,3,5) \\
12 && (1,2,3) && 1,1,-1,-1 & (1,2,3) & \text{itself} \\
13 && (1,3,4) &1,1,-1,-2&& \tfrac{1}5(1,2,3) & \PP(1,3,4,5) \\
14 && (1,2,5) &1,1,-1,-3&& (1,1,2) & \PP(1,2,3,5) \\
15 && (1,2,7) &1,1,-1,-5&& \tfrac{1}3(1,1,2) & \PP(2,3,5,7) \\
 &\multicolumn{5}{l|}{\quad\text{and the inverse of }3} & \\
\cmidrule(lr){1-7}
16 & \PP(1,1,2,3) & \tfrac{1}3(1,1,2) &&&\text{Mfs} & \PP^1 / \PP(1,1,2) \\
17 && \tfrac{1}2(1,1,1) && 3,1,-1,-1 &\text{Mfs} & \PP^2 / \PP^1 \\
18 && (1,1,2) & 3,1,-1,-5&& \tfrac{1}3(1,1,2) & \PP(1,2,3,5) \\
19 && (1,1,3) & 2,1,-1,-3 &&\text{Mfs} & \PP^2 / \PP^1 \\
20 && (1,2,3) &&& \text{Mfs} & \PP^1 / \PP(1,2,3) \\
21 && (1,2,3) & 3,2,-1,-5&& (1,1,1) & \PP(1,2,3,5) \\
22 && (1,1,4) & 2,1,-1,-5&& (1,1,1) & \PP(1,3,4,5) \\
23 && (1,1,5) & 2,1,-1,-7&& \tfrac{1}2(1,1,1) & \PP(2,3,5,7) \\
24 && (1,3,4) &&1,1,-1,-1& (1,3,4) & \text{itself} \\
25 && (1,3,5) & 1,1,-1,-2&& (1,2,3) & \PP(1,2,3,5) \\
26 && (1,4,5) & 1,1,-1,-3&& \tfrac{1}7(1,3,4) & \PP(2,3,5,7) \\
27 && (1,3,7) & 1,1,-1,-4&& \tfrac{1}5(1,2,3) & \PP(3,4,5,7) \\
 &\multicolumn{5}{l}{\quad\text{and the inverses of }4, 8, 10} & \\
\cmidrule(lr){1-7}
28 & \PP(1,2,3,5) & \tfrac{1}5(1,2,3) &&& \text{Mfs} & \PP^1 / \PP(1,2,3) \\
29 && (1,1,3) & 3,1,-1,-4 && (1,1,2) & \PP(1,3,4,5) \\
30 && (1,1,4) & 3,1,-1,-7 && \tfrac{1}3(1,1,2) & \PP(3,4,5,7) \\
31 && (1,2,5) & 2,1,-1,-5 && \tfrac{1}5(1,2,3) & \PP^3 / \tfrac{1}5(1,2,3,4) \\
 &\multicolumn{5}{l|}{\quad\text{and the inverses of }11, 14, 18, 21, 25} & \\
\cmidrule(lr){1-7}
32 & \PP(1,3,4,5) & \tfrac{1}4(1,1,3) && 5,1,-1,-3 & \text{Mfs} & \PP^2 / \PP^1 \\
33 && (1,1,2) & 5,1,-1,-7 && \tfrac{1}5(1,1,4) & \PP(2,3,5,7) \\
34 && (1,1,3) & 4,1,-1,-7 && \tfrac{1}4(1,1,3) & \PP(3,4,5,7) \\
35 && (1,2,3) & 3,1,-1,-5 && \tfrac{1}7(1,2,5) & \PP(3,4,5,7) \\
 &\multicolumn{5}{l|}{\quad\text{and the inverses of }5, 13, 22, 29} & \\
\cmidrule(lr){1-7}
 &\multicolumn{4}{l}{\PP(2,3,5,7)\quad\text{the inverses of }15, 23, 26, 33} & \\
\cmidrule(lr){1-7}
 & \multicolumn{4}{l}{\PP(3,4,5,7) \quad\text{the inverses of }27, 30, 34, 35} & \\
\cmidrule(lr){1-7}
 & \multicolumn{4}{l}{\PP^3 / \tfrac{1}5(1,2,3,4)\quad\text{the inverse of }31} & \\
\bottomrule
\end{array}
\]
\caption{Sarkisov links between the~$8$ toric Mori\nobreakdash--Fano~$3$-folds and some Mori fibre spaces~(Mfs). We abbreviate a divisorial contraction~$(-r,a,b,c)$ by~$\tfrac{1}r(a,b,c)$, omitting the fraction when~$r=1$. Mfs are indicated only by~$F/B$, where~$F$ is the fibre and~$B$ is the base.}
\label{fig!links}
\end{table}
\subsection{The midpoints of Sarkisov links}\label{sec!midpoints}
The endpoints of a Sarkisov link are not the only places a Fano variety might appear. The Sarkisov link~\eqref{eq!link} may include a flop~$Y_i\rightarrow Z_i \leftarrow Y_{i+1}$. In that case,~$Z_i$ is a toric Fano \nf\ with terminal singularities and Picard rank~$\rho_{Z_i}=1$, but with Weil divisor class group of rank~2 so not \Qfact. Karzhemanov's famous example~$Z_{70}\subset\PP^{37}$~\cite[1.2]{karzhemanov64}, the unique Fano \tf\ of degree~$-K_Z^3=70$, is of this nature, though it is not toric. Starting with~$X=\PP(1,1,4,6)$, Karzhemanov makes a (weighted) blowup~$X\leftarrow Y_1$ of a midpoint of the~$\PP(4,6)$ stratum~$L\subset X$ to give a weak Fano \tf\ $Y_1$ on which the proper transform of~$L$ is a flopping curve:~$Z_{70}$ is the anticanonical image of~$Y_1$, namely the base of that flop.

The links we construct for toric Fano \ffs\ provide around~$10,\!000$such flopping bases, each of which is a toric Fano \ff, the majority not known to us. Of these, around~$2000$ arise as~$Z_1$ in~\eqref{eq!link}, that is from a terminal weak Fano blowup of a (possibly fake) weighted projective space. For example, if~$X=\PP(1^4,2)$ with orbifold point~$Q=(0:0:0:0:1)\in X$, and~$X\leftarrow Y_1$ is the blowup of a smooth point~$P\in X$, then the proper transform~$C\subset Y_1$ of the line~$L\cong\PP(1,2)$ through~$P$ and~$Q$ is a flopping curve (and a toric stratum, for suitable torus). After contracting~$C$, the result is a Gorenstein Fano \ff\ $Z_1\subset\PP^{114}$ of degree~$567$ in its anticanonical embedding.

It is also possible for one of the \Qfact, Picard rank~2 varieties~$Y_i$ appearing in a link~\eqref{eq!link} to be a Fano variety. Indeed if the link starts with a flip~$Y_1\dashrightarrow Y_2$, then~$Y_1$ is such a Fano variety, or if there is a sequence~$Y_{i-1}\dashrightarrow Y_i\dashrightarrow Y_{i+1}$ comprising an antiflip followed by a flip, then again~$Y_i$ is a Fano variety; we give an example of this in~\S\ref{sec!12345} below.
\section{Operations on a simplicial fan}\label{sec!operations}
We summarise the toric operations that arise in the links we construct, closely following Reid \cite[\S\S2--4]{reid82}, which also has guiding pictures. We start with a terminal Fano (fake) weighted projective~$n$\nobreakdash-space~$X$ determined by a fan~$\Si$ of cones~$C_1,\dots,C_{n+1}$ on rays generated by vertices~$v_1,\dots,v_{n+1}$ that are primitive (indivisible in~$L$). In the notation of~\eqref{eq!link}, so that, for example, each~$Y_i$ is \Qfact\ of Picard rank~$\rho_{Y_i}=2$, the steps of the link are:

\renewcommand{\labelenumi}{(\arabic{section}.\arabic{enumi})}
\begin{enumerate}
\item Terminal extractions.\label{type:extraction}
\end{enumerate}
The map~$X\leftarrow Y_1$ arises from a fan refinement~$\Si_1\subset\Si$, where~$\Si_1$ is the subdivision of~$\Si$ by the ray though a new primitive vertex~$v\in L$.

In dimension~$3$, if~$v\in C_i$ and~$C_i$ is a terminal quotient singularity~$\tfrac{1}{r_i}(1,a,r_i-a)$ of index~$r_i>1$, then~$v$ is necessarily the vertex of height~$r_i+1$ inside~$C_i$, the so-called~\emph{Kawamata blowup}~\cite{kawamata96} (denoted by~$\tfrac{1}{r_i}(1,a,r_i-a)$ in Table~\ref{fig!links}). If~$C_i$ is a nonsingular cone (its vertices~$s_1,s_2,s_3$ form a basis of~$L$) then either~$v = s_1 + as_2 + bs_3\in C_i^\circ$ for coprime~$a,b>1$ (denoted by~$(1,a,b)$ in Table~\ref{fig!links}) or~$v=s_1+s_2\in\partial C_i$, up to permutations of the~$s_i$ (denoted by~$(1,1,0)$ in Table~\ref{fig!links}). In dimension~$4$ there are no general results specifying those~$v\in C_i$ that determine terminal extractions.

\begin{enumerate}[resume]
\item Flips, flops and anti-flips~\cite[\S3]{reid82}.\label{type:flip}
\end{enumerate}
A map~$\phi_i\colon Y_i\rightarrow Z_i$ may arise by amalgamating a union of~$\ge 2$ cones~$C_j$ of~$\Si_i$ into a single convex cone~$D$. If~$D$ is not simplicial, then in our context~\eqref{eq!link}~$D$ necessarily has~$n+1$ vertices and (as~$\rho_{Y_i}=2$ so~$\rho_{Y_i/Z_i}=1$) there is a unique alternative way of subdividing~$D$ into simplicial cones on the same vertices. This forms a new simplicial fan~$\Si_{i+1}$, giving the composition~$Y_i\rightarrow Z_i\leftarrow Y_{i+1}$ which is an isomorphism in codimension~$1$ (see~\cite[(3.4)]{reid82}).

Consider the vertices~$\{s_1,\dots,s_{n+1}\}$ of~$D$, ordered so that~$s_1,\dots,s_n$ span one of the cones~$C_j$ of~$\Si_i$ and thus lie on an affine hyperplane~$H\subset L$. If~$s_{n+1}$ also lies on~$H$, then~$Y_i\dashrightarrow Y_{i+1}$ is a flop ($\phi_i$~contracts only curves~$\Ga$ with~$K\Ga=0$) and~$Y_{i+1}$ again has \Qfact\ terminal singularities. If~$s_{n+1}$ lies on the same side of~$H$ as the origin, then~$-K_{Y_i}$ is~$\phi_i$\nobreakdash-ample,~$Y_i\dashrightarrow Y_{i+1}$ is a flip, and again~$Y_{i+1}$ has \Qfact\ terminal singularities (see~\cite[\S4]{reid82}).

On the other hand, if~$s_{n+1}$ lies on the opposite side of~$H$ to the origin, then~$Y_i\dashrightarrow Y_{i+1}$ is an anti-flip, and we lose control of the singularities of~$Y_{i+1}$: the shed has grown and may now admit interior points. We must check: if the singularities are terminal, then we continue with~$Y_{i+1}$; if not, then the link is bad.

In Table~\ref{fig!links} each of these operations is denoted by a vector~$(b_1,\dots,b_{n+1})$ for which~$\sum b_ks_k=0$.

\begin{enumerate}[resume]
\item Divisorial contractions. \label{type:divisor}
\end{enumerate}
Continuing the notation of (\ref{sec!operations}.\ref{type:flip}), if~$D$ is simplicial then~$\phi_i$ is a divisorial contraction to a point and~$Z_i=X'$ is the end of the link (again as~$\rho_{Y_i}=2$ so~$\rho_{Y_i/Z_i}=1$).

It may also happen that~$\phi_i\colon Y_i\rightarrow Z_i$ arises by a subset of cones of~$\Si_i$ combining to make a union of simplicial cones of a new fan, in which case~$\phi_i$ is a divisorial contraction to a locus of dimension~$>0$, and again we have reached the end of the link.

\begin{enumerate}[resume]
\item Mori fibre spaces~\cite[(2.5--6)]{reid82}. \label{type:mfs}
\end{enumerate}
When~$\dim Z_i<\dim Y_i$, then~$i=s$ is the end of the link and~$Y_s\rightarrow Z_s$ is a Mori fibre space with (possibly fake) weighted projective spaces as fibres.

\begin{enumerate}[resume]
\item{Blowup and flip notation and the weights of~$\C^*$ variations.}
\end{enumerate}
As indicated above, Table~\ref{fig!links} abbreviates the data of each birational map by a sequence of numbers, namely the coefficients of the linear relations among rays of the relevant cones. For \ff\ output, including the much larger data set recorded at~\cite{grdbweb}, we abbreviate the birational maps as follows.

Blowups (and blowdowns) are determined by their centre~$P$ and the primitive lattice point~$v$ on the the new subdividing ray that describes the blown-up fan. Thus we record blowups (and blowdowns) by two pieces of data, as follows: first, a sequence~$[a_0,a_1,\dots,a_d]$ indicating the weights of coordinates along the blowup centre -- that is, the~$d$\nobreakdash-dimensional toric stratum~$P=\PP(a_0,\dots,a_d)$ being blown up -- and, second, a vector~$(b_0,b_1,\dots,b_{4-d})$ that gives the coefficients of the minimal integral relation defining~$v$ in terms of the rays of~$P$ (in some order that is not recorded). The notation does not determine the (singularity) type of the centre being blown up, other than its dimension and lattice index, but that can be recovered unambiguously from the weighted projective space.

Some examples illustrate this notation for blowups of a toric \ff:
\begin{itemize}
\item
$[1](-1,1,1,1,2)$ is the~$(1,1,1,2)$-weighted blowup of a nonsingular toric point-strat\-um. The blowup is determined by the fan subdivision at a ray through~$v = s_1+s_2+s_3+2s_4\in L$, where~$s_i$ are the vertices of a regular~$4$\nobreakdash-dimensional cone corresponding to the point, while the sequence~$[1]$ indicates that this is a cone of index~$1$, that is, a regular cone.
\item
$[1,1](-1,1,1,1)$ is the ordinary blowup of a nonsingular toric curve stratum The blowup is determined by the subdivision at~$v=s_1+s_2+s_3$, which lies on a regular $3$\nobreakdash-dimensional cone with vertices~$s_i$.
\item
$[3](-3,1,1,1,2)$ is the blowup by~$v=\tfrac{1}3(s_1+s_2+s_3+2s_4)$ in a~$4$\nobreakdash-dimensional cone of index~$3$.
\item
$[4](-2,1,1,1,1)$ is the blowup by~$v=\tfrac{1}2(s_1+\cdots+s_4)$ in a~$4$\nobreakdash-dimensional cone of index~$4$, for which~$\sum s_i$ is divisible by~$2$ in the lattice~$L$, but not by~$4$.
\item
$[2,2](-2,1,1,1)$ is the Kawamata blowup along a curve of transverse~$\tfrac{1}2(1,1,1)$ singularities, that is, the blowup by~$v=\tfrac{1}2(s_1+s_2+s_3)$ in the corresponding~$3$\nobreakdash-dimensional cone.
\item
$[2,4](-2,1,1,1)$ is the unique extremal extraction along a curve of generically transverse~$\tfrac{1}2(1,1,1)$ singularities that equals the Kawamata blowup at a general point. (This occurs, for example, as the start of a link from~$\PP(1,1,2,3,4)$.)
\end{itemize}
For flips and other isomorphisms in codimension~1, we record a sequence of integers that are the coefficients of the minimal integral relation among the rays of~$D$, in the notation of~(\ref{sec!operations}.\ref{type:flip}). They may also be treated as weights for a~$\C^*$ action; see~\cite[\S1]{B99} for example.
\begin{itemize}
\item
$(4,1,-1,-1,-3)$ is a \ff\ flop contracting~$\PP(4,1)$ and~$\PP(1,1,3)$ on the two sides respectively to a common point in the base.
\item
$(2,1,0,-1,-1)$ is a \ff\ flip over a curve in the base that makes a Francia flip (namely $(2,1,-1,-1)$ in a~$3$-fold) over each point of the curve.
\end{itemize}
In each case, the rays involved may or may not generate the whole lattice, but we do not record this co-index. 

\section{Extending Prokhorov's web}\label{sec!dim3}
Prokhorov~\cite[\S10]{prokhorov10} computes links from \tf\ weighted projective spaces that start with the Kawamata blowups of quotient singularities. We list all toric Sarkisov links between toric Mori\nobreakdash--Fano \tfs\ in Table~\ref{fig!links}, and discuss some particular cases here. Since we list the \tfs\ $X$ with smaller weights first (as one might if using Sarkisov links to construct `complicated' Fano \tfs\ from `simple' ones), the typical behaviour is a relatively high disrepancy blowup followed by an antiflip; the picture would be reversed if we listed larger weights first. By~\cite[5]{kawamata96}, it remains only to consider the blowups of~$1$\nobreakdash-strata that do not pass through a quotient singularity and weighted blowups~$(1,a,b)$ of a smooth~$0$\nobreakdash-stratum. The latter is an infinite collection of blowups, and we discuss bounds in~\ref{sec!bounds}.

\subsection{Links from a smooth point.}
Consider the weighted~$(1,1,2)$ blowup of the smooth~$0$\nobreakdash-stratum of~$\PP(1,3,4,5)$. This extends to a Sarkisov link:
\[
\begin{array}{cccccccccc}
 \multicolumn{2}{r}{\buildrel{(-1,1,1,2)}\over{}}& Y_1 &&\buildrel{(4,1,-1,-3)}\over{\dashrightarrow}&& Y_2&\multicolumn{3}{l}{\buildrel{(3,1,1,-1)}\over{}}\\
 &\swarrow&&\searrow && \swarrow&&\searrow \\
\PP(1,3,4,5) &&&& Z_1 &&&& \PP(1,2,3,5)
\end{array}
\]
In terms of the rays of the fan, we start with~$\{(1,1,0),(0,-1,1),(1,-1,-1),(-2,0,-1)\}$, and the blowup inserts
\[
(0,1,0) = 1\cdot (-2,0,-1) + 1\cdot (0,-1,1) + 2\cdot (1,1,0). \tag{$-1,1,1,2$}
\]
The flip then expresses the two subdivisions
\[
4\cdot (0,1,0) + 1\cdot (1,-1,-1) = 1\cdot (-2,0,-1) + 3\cdot (1,1,0) \tag{$4,1,-1,-3$}
\]
after which the ray~$(1,1,0)$ can be contracted.
\subsection{Links from the eighth toric Mori\nobreakdash--Fano 3-fold.}
We denote the action of~$\eps\in \Z/5$ act on~$\PP^3$ by~$(\eps,\eps^2,\eps^3,\eps^4)$ by~$\tfrac{1}5(1,2,3,4)$. Then~$\PP^3 / \tfrac{1}5(1,2,3,4)$ is a well-known toric Mori\nobreakdash--Fano \tf\ which is a fake weighted projective space (compare~\cite{bb92}). Its only toric extremal extraction is the~$\tfrac{1}5(1,2,3)$ blowup of any of the four quotient singularities, and this extends to a Sarkisov link:
\begin{equation}\label{eq!P35}
\begin{array}{cccccccccc}
 \multicolumn{2}{r}{\buildrel{(-5,1,2,3)}\over{}}& Y_1&\buildrel{(5,1,-1,-2)}\over{\dashrightarrow}& Y_2&\multicolumn{3}{l}{\buildrel{(5,2,1,-1)}\over{}}\\
 &\swarrow&&&&\searrow \\
X = \PP^3 / \tfrac{1}5(1,2,3,4) &&&&&& \PP(1,2,3,5) = X'
\end{array}
\end{equation}
In coordinates,~$X$ may be defined by the four cones on the vertices
\[
v_1=(0,1,1),\quad
v_2=(-1,0,-2),\quad
v_3=(-1,-2,1),\quad
v_4=(2,1,0)
\]
(whose sum is the origin, but which generate only a sublattice of index~$5$),~$Y_1$ is the~$\tfrac{1}5(1,2,3)$ blowup at the new vertex~$v_5=(1,1,0)$ which satisfies
\[
5\cdot v_5 =1\cdot v_2 + 2\cdot v_1+3\cdot v_4 \tag{$-5,1,2,3$}
\]
$Y_2$ is the flip by subdividing~$\left<v_1,v_3,v_4,v_5\right>$ the other way, corresponding to the linear relation
\[
5\cdot v_5 + 1\cdot v_3 = 1\cdot v_1 + 2\cdot v_4 \tag{$5,1,-1,-2$}
\]
and~$X'=\PP(1,2,3,5)$ follows from blowing down
\[
5\cdot v_5 + 2\cdot v_3 + 1\cdot v_2 = v_4 \tag{$5,2,1,-1$}
\]
leaving
\[
1\cdot v_1 + 2\cdot v_2 + 3\cdot v_3 + 5\cdot v_5 = 0. \tag{$1,2,3,5$}
\]
\begin{figure}[ht]
\includegraphics[scale=0.53]{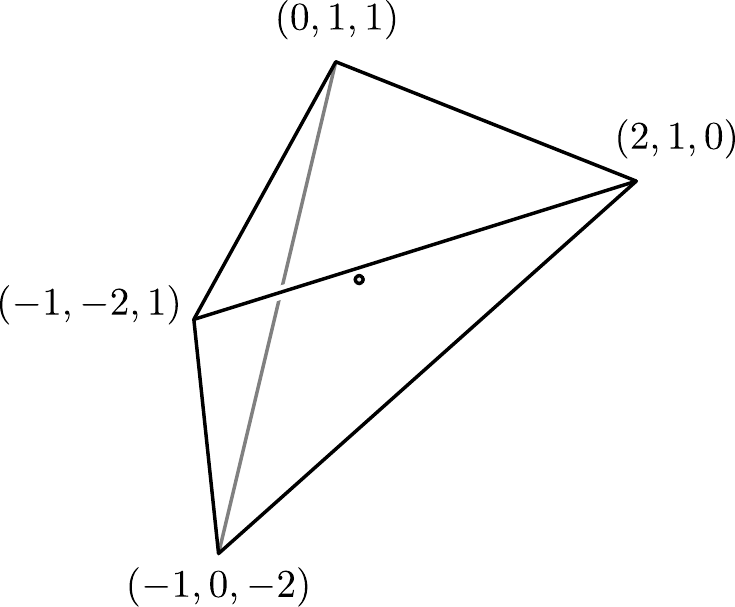}
\includegraphics[scale=0.53]{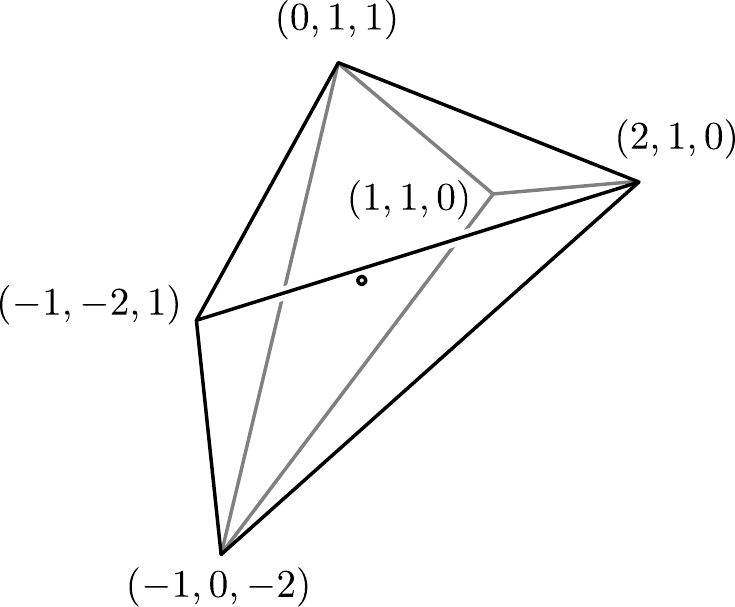}
\includegraphics[scale=0.53]{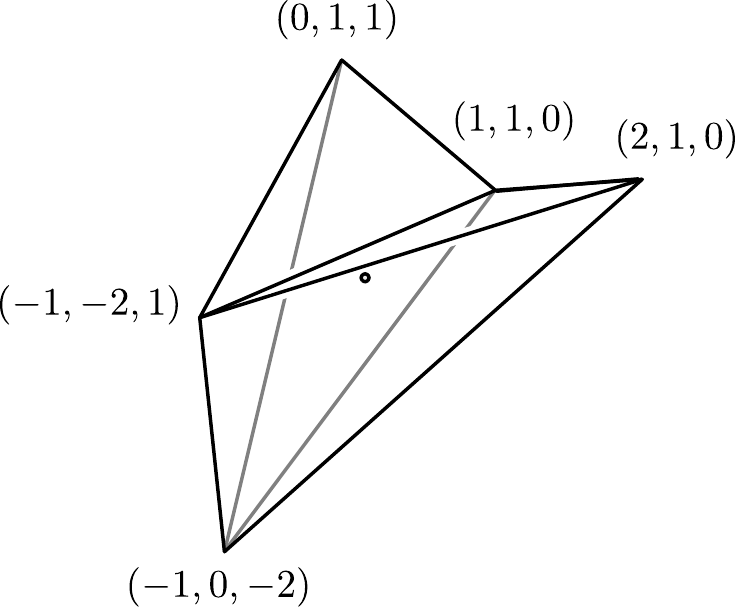}
\includegraphics[scale=0.53]{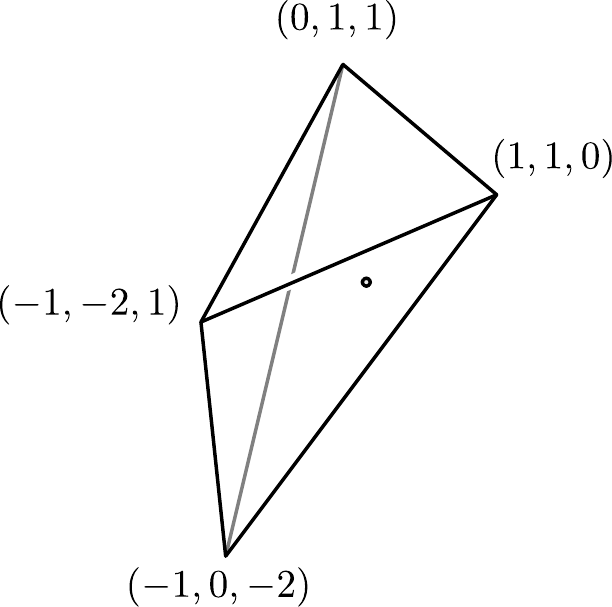}
\caption{The sheds of the toric fans of~$X\leftarrow Y_1\dashrightarrow Y_2\rightarrow X'$ in the link~\eqref{eq!P35}\label{fig!pic}}
\end{figure}
\subsection{Links from 1-dimensional centres.}
A toric~$1$\nobreakdash-stratum~$L\subset X$ is the centre of an extremal extraction only when~$X$ is smooth in a neighbourhood of~$L$ and~$X\leftarrow Y_1$ is the blowup of~$L$. In each case there is only a single~$L\subset X$ up to symmetry, and its blowup extends to a link; the blowup (or blowdown) is denoted~$(1,1,0)$ in Table~\ref{fig!links}.
\subsection{Bounding links from a smooth centre.}\label{sec!bounds}
Though the number of Sarkisov links between toric Fano \tfs\ is finite, the number of terminal extractions from a smooth point is infinite, even in the toric case: the~$(1,a,b)$ blowup is terminal for coprime~$a,b\ge1$.

Consider the main case of the blowup of a point on~$\PP^3$. We may define the fan of~$\PP^3$ on the four rays~$e_1,e_2,e_3$ (the standard basis of~$L$) and~$e_4=(-1,-1,-1)$, and without loss of generality we consider blowing up by the point~$e_5=(a,1,b)$ with~$a>b>1$ coprime. (By symmetry, the only other cases are when~$b=1$, and these smaller cases may be handled similarly.) Firstly, the blowup is indeed terminal. After blowup, the edge~$e_1,e_3$ must be antiflipped to~$e_4,e_5$: indeed the equation
\[
g_{134} = x - 3y + z
\]
which supports the roof of the shed~$(g_{134}=1)$ of cone~$\si_{134}=\left< e_1,e_3,e_4\right>$ is strictly positive at~$e_5$, as~$a\ge3$ and~$b\ge2$. This shows that~$-K$ is positive on the toric~$1$\nobreakdash-stratum corresponding to~$\left<e_1,e_3\right>$, but even without that observation the union of cones~$\si_{134}\cup\si_{135}$ (which is strictly convex, as~$g_{134}>0$) may be differently subdivided into convex cones by the hyperplane through~$e_4$ and~$e_5$.

After this antiflip, the lattice point~$v=(1,0,1)$ lies in the cone~$\si_{345}=\left< e_3,e_4,e_5 \right>$. This cone has supporting equation~$(g_{345}=1)$ where
\[
g_{345} = \left(\frac{3-b}{a-1}\right)x + \left(\frac{b-2a-1}{a-1}\right) y + z
\]
so that, if~$a>b\ge3$, the blow up at~$v$ has discrepancy~$g_{345}(v)-1=(3-b)/(a-1)\le0$, which violates terminality.

In the case~$a>b=2$, the link would end by contracting
\[
e_1 = \tfrac{1}{a-2}(1\cdot e_2 + 2\cdot e_4 + 1\cdot e_5)
\]
but that contraction to a quotient singularity~$\tfrac{1}{a-2}(1,2,1)$ has non-positive discrepancy unless~$a-2\le3$. Thus we must be in the situation~$a\le5$ and~$b=2$, and the possible blowups leading to a Sarkisov link are indeed bounded.
\subsection{Higher rank Fano 3-folds, links and relations.}
Since each link~$X\dashrightarrow X'$ in Table~\ref{fig!links} identifies the big torus, it is easy to compose sequences of these links to give birational automorphims~$X\dashrightarrow X$, describing relations in the Sarkisov program. Following~\cite{kaloghiros13}, relations are derived from minimal model programs~(MMPs) on certain varieties~$Z$ with~$\rho_Z=3$. Explicitly in the toric case, most links in Table~\ref{fig!links} arise by patching together two different~MMPs on toric Fano \tfs\ $Y$ with~$\rho_Y=2$ (since one of the~$Y_i$ in~\eqref{eq!link} is Fano, unless there is a flop). Such~$Y$ are classified (into~$35$ cases, coincidentally), and further~MMPs from the~$75$ rank~$3$ toric Fano \tfs\ describe relations. This can all be carried out by similar, but bigger, toric considerations.
\section{The web of 4-dimensional toric Sarkisov links}\label{sec!dim4}
The main result of this paper is an attempt to construct data analogous to Table~\ref{fig!links} for toric Mori\nobreakdash--Fano \ffs. We do not attempt a complete classification of links; compare~\cite{tiago} for the complete classification in the case~$X=\PP^4$. Instead, for each of the~$35,\!947$ varieties we compute all toric (point, curve and surface) blowups of discrepancy at most~$5$, and for each one either extend to a Sarkisov link or discard as a bad link. Even subject to this severe constraint on discrepancy (at one end of the link) there are over a million links, and we record them all in full detail online at~\cite{grdbweb}. We illustrate some general features of the results here by discussing links from three exemplary starting points~$X$.
\subsection{Blowups of~$\PP(1,2,3,4,5)$}\label{sec!12345}
Under the given restrictions on blowups and their discrepancies, there are~$275$ Sarkisov links from~$X=\PP(1,2,3,4,5)$.

For example, there are~$4$ ways of making a link from a~$(-5,1,2,3,4)$ blowup of the index~$5$ point: after the blowup~$X\leftarrow Y_1$, these links are completed by
\begin{enumerate}
\renewcommand{\labelenumi}{(\roman{enumi})}
\item
$Y_1\rightarrow \PP(1,2,3,4)$ Mori fibre space with~$\PP^1$ fibre
\item
a~$(3,1,-1,-1,-2)$ flop, then~$(4,3,1,-1,-2)$ flip to~$Y_3\rightarrow \PP^1$ Mori fibre space with $\PP(1,1,1,2)$ fibre
\item
a~$(2,1,0,-1,-1)$ flip, then~$(1,2,3,4)$-weighted blowdown~$Y_2\rightarrow\PP(1,1,1,2,3)$ to a smooth point
\item
a~$(4,1,-1,-2,-3)$ antiflip, then~$(3,2,1,-1,-2)$ flip, then~$(1,2,2,3)$-weighted blowdown $Y_3\rightarrow\PP(1,1,2,3,4)$ to a smooth point.
\end{enumerate}
The last of these exhibits antiflip--flip behaviour which are rare for \tfs\ (the index~$1$ cases in~\cite{CM,BZ,Campo} all involve a flop, but Guerreiro has examples for special members of a family in higher index). In particular, the variety~$Y_2$ in the middle of the link is a Fano \ff: it has Picard rank~$2$, with both extremal rays flipping, and may be defined as the~$(\C^\times)^2$ quotient
\[
\left(
\begin{array}{ccc|ccc}
0&1&2&3&4&1\\1&2&3&4&5&0
\end{array}
\right)
\]
This phenomenon seems to be common in dimension~$4$.
\subsection{Links between fake weighted projective spaces}
The web of Sarkisov links is of course connected, and so there are necessarily links between fake weighted projective spaces with different discriminant groups. For example there is a Sarkisov link
\[
X=\PP(2,3,5,5,13)/\tfrac{1}5(0,1,3,4,3)
\dashrightarrow
X'=\PP(4,5,6,7,17) / \tfrac{1}2(1,1,1,0,0)
\]
that factorises as
\[
\begin{array}{ccccccccccccccccc}
 \multicolumn{2}{r}{\buildrel{(-25,1,5,6,14)}\over{}}& Y_1 &&\buildrel{(15,1,-1,-2,-8)}\over{\dashrightarrow}&& Y_2 && \buildrel{(65,6,1,-7,-34)}\over{\dashrightarrow} && Y_3 & \multicolumn{2}{l}{\buildrel{(25,7,5,2,-12)}\over{}}\\
&\swarrow&&&&&&&&&&\searrow & \\
X &&&&&&&&&&&& X'
\end{array}
\]
that is, a divisorial extraction from a point, two consecutive flips, followed by a divisorial contraction to a point. Again, this phenomenon is not something we have seen in~$3$ dimensions (Campo~\cite{Campo} finds cases with multiple flips, but always following a flop). In particular, the variety~$Y_1$ is a Fano \ff\ with Picard rank~$2$.
\subsection{Blowups of~$\PP^4$}\label{P4blowups}
Starting with~$X=\PP^4$, our bounds permit just one example in which the initial extremal extraction is followed by a flop
\[
\begin{array}{cccccccccccc}
 \multicolumn{2}{r}{\buildrel{(-1,1,1,2,2)}\over{}}& Y_1&\buildrel{(1,1,0,-1,-1)}\over{\dashrightarrow}& Y_2&\multicolumn{3}{l}{\text{\small ($\PP(1,1,1,2)$-bundle)}}\\
 &\swarrow&&&&\searrow \\
\PP^4 &&&&&& \PP^1
\end{array}
\]
However, allowing bigger discrepancies it is easy to construct cases of the form
\[
\PP^4 \leftarrow Y_1 \buildrel{\text{flop}}\over{\dashrightarrow} Y_2
\buildrel{\text{flip}}\over{\dashrightarrow} Y_3 \rightarrow X'
\]
where the blowups have weights
\begin{equation}\label{eq!abcd}
(7, 8, 9, 12),\quad (10, 11, 13, 17),\quad (11, 12, 14, 19),\quad (13, 14, 17, 22)
\end{equation}
and the resulting~$X'$ are respectively
\[
\PP(1,3,4,5,12),\quad \PP(1,4,6,7,17),\quad \PP(1,5,7,8,19),\quad \PP(1,5,8,9,22).
\]
(The~$4$-tuples in~\eqref{eq!abcd} are simply the solutions~$(d,a,b,c)$ of the equation~$a+b+c=4d+1$ with~$a,b,c\ge d\ge1$, with any three coprime, for which the blowup~$(a,b,c,d)$ of the standard toric~$\PP^4$ terminal, subject to the bound~$a,b,c,d\le100$. The equation guarantees the flop.)

Similarly, our bounds permit only a single case when the blowup is followed by a flip
\[
\begin{array}{cccccccccccc}
 \multicolumn{2}{r}{\buildrel{(-1,3,3,4,5)}\over{}}& Y_1&\buildrel{(3,1,0,-1,-2)}\over{\dashrightarrow}& Y_2&\multicolumn{2}{l}{\buildrel{(4,1,1,1,-1)}\over{}}\\
 &\swarrow&&&&\searrow \\
\PP^4 &&&&&& \PP(1,1,2,2,5)
\end{array}
\]
Again higher discrepancies allow a handful of other cases~$\PP^4\leftarrow Y_1\buildrel{\text{flip}}\over{\dashrightarrow} Y_2\rightarrow X'$, namely weighted blowups~$(4,4,5,7)$,~$(5,5,6,8)$,~$(12,13,15,20)$ with~$X'$ respectively~$\PP(1,2,3,3,7)$, $\PP(1,2,3,3,8)$ and $\PP(1,5,7,8,20)$. (These are the solutions~$(d,a,b,c)$ to~$a+b+c<4d+1$ with~$a,b,c\le100$,~$d\le50$ and the same additional conditions as above.) In each case,~$Y_1$ is a Fano \ff\ of Picard rank~$2$.

The typical behaviour, though, is that the extremal extraction is followed by an antiflip. This can be followed by a flop or a flip, but a sequence of two antiflips is common, as in the following:
\[
\begin{array}{ccccccccccccccccc}
 \multicolumn{2}{r}{\buildrel{(-1,1,2,5,9)}\over{}}& Y_1 &&\buildrel{(1,1,-1,-4,-8)}\over{\dashrightarrow}&& Y_2 && \buildrel{(2,1,-1,-3,-7)}\over{\dashrightarrow} && Y_3 & \multicolumn{2}{l}{\buildrel{(5,4,3,1,-4)}\over{}}\\
&\swarrow&&&&&&&&&&\searrow & \\
\PP^4 &&&&&&&&&&&& \PP(1,4,7,8,9)
\end{array}
\]
In this case,~$Y_3$ is a Fano \ff\ of Picard rank~$2$.
\section{Further related problems}
\subsection{High-discrepancy blowups.}\label{sec!highdiscrep}
The list we describe above certainly misses some Sarkisov links, such as those in~\S\ref{P4blowups}, since we restrict attention to discrepency at most~$5$ at one of the ends. There are only finitely many Sarkisov links from any given centre, and with more work one can bound the terminal extractions that initiate these links and describe all the missing cases. In the case of~$X=\PP^4$, Guerreiro~\cite{tiago} computes all~$421$ Sarkisov links from~$\PP^4$ that start with the blowup of a point.
\subsection{Non-terminal singularities.}
The classification~\cite{kasprzyk10} includes a further~$348,\!930$ Fano \tfs\ $X$ with strictly canonical singularities and~$\rho_X=1$. It makes sense to compute links of the form~\eqref{eq!link} between these, for example by permitting crepant extractions, even though they would be regarded as `bad' links in the usual terminal situation.
\subsection{Running the Sarkisov program on toric~$n$-folds.}
In the toric context, any two toric Mori\nobreakdash--Fano \nfs\ $X$ and~$Y$ have a common big torus, and any such identification extends to a birational map~$X\dashrightarrow Y$. At first sight, given any two Fano \nf\ polytopes~$P_1$ and~$P_2$ and elements~$g_1,g_2\in\GL(n,\Z)$, the task is to compute a series of operations of the form~(\ref{sec!operations}.\ref{type:extraction}--\ref{type:divisor}) that takes the spanning fan of~$g_1(P_1)$ to that of~$g_2(P_2)$, via fans on at most~$n+1$ rays. But the Picard rank can increase when Mori fibre spaces appear, and the analysis may be rather subtle; the Sarkisov program applies to Mori fibre spaces, not only Mori\nobreakdash--Fano \nfs. This is a much more substantial problem than the one we address in~\S\ref{sec!dim4}.
\subsection{Higher Picard rank.}
The secondary fan of a toric variety~$V$ (of arbitrary Picard rank) contains the data of all minimal model programs~(MMPs). Magma~\cite{magma} includes functions to run all toric~MMPs from a given~$V$, inductively contracting all extreme rays of the Mori cones that arise. The sets of all toric varieties and all their~MMPs are infinite but could be enumerated up to a bound. This would contain all Sarkisov links that are dominated by such~$V$ and, following~\cite{kaloghiros13}, relations in the toric Sarkisov program among those.
\subsection{Fano 3-folds and 4-folds more generally.}
The toric case considered here is a test case for overlaying other (partial) classifications of Fano varieties by webs of Sarkisov links. Analysis of the numerical data of possible links has had spectacular success contributing to the classification of particular types of Fano varieties (for example, in works of Takeuchi, Alexeev, Takagi and Prokhorov, among others) and a `numerical' web in the spirit of the Graded Ring Database~\cite{grdbweb} may be possible.
\subsection*{Acknowledgements}
This work was initiated as part of EPSRC grant EP/E000258/1 and completed with the support of EPSRC fellowship EP/N022513/1. Buczy{\'n}ski was partially supported by the National Science Center, Poland, project 2017/26/E/ST1/00231. The computations were performed using the Imperial College High Performance Computing Service and the Compute Cluster at the Department of Mathematics, Imperial College London, and we thank Simon Burbidge, Matt Harvey and Andy Thomas for valuable technical assistance. We thank John Cannon and the Computational Algebra Group at the University of Sydney for providing licenses for the computer algebra system Magma, and for support and hospitality.
\end{document}